\documentclass{article}
\usepackage{graphicx}
\usepackage{amssymb,amsmath}
\def\Aut{\operatorname{Aut}}
\def\ds{\displaystyle}

\def\timehm{\count31=\time \count32=\count31 \divide\count31 by 60
\number\count31 \multiply\count31 by 60 \advance\count32 by
-\count31 :\ifnum\count32<10 0\fi \number\count32}

\def\ep{\varepsilon}

\def\ap{\alpha}

\def\be{\beta}

\def\tm#1{\left(\text{mod }#1\right)}

\newtheorem{thm}{Theorem}

\newtheorem{lem}[thm]{Lemma}

\newtheorem{conj}[thm]{Conjecture}

\newcommand{\qed}{\hfill\rule{2mm}{2mm}}

\newcommand{\leg}[2]{ \left( \frac{#1}{#2} \right) }

\def\sod{^{\kern 14pt odd}}
\def\be{\mathfrak b}

\def\End{\mathrm{End}}
\def\Gal{\mathrm{Gal}}
\def\GL{\mathrm{GL}}
\def\mod{\;{\rm mod}\;}

\def\tr{\mathop{\rm tr}\nolimits}

\def\Q{{\cal Q}}
\def\F{\mathbb{F}}

\def\Q{\mathbb{Q}}

\def\Z{\mathbb{Z}}

\def\ideal#1{<\kern-2pt #1\kern-2pt >}

\title{Square-free discriminants of Frobenius rings}
\begin{author}{Chantal David and
Jorge Jim\'enez Urroz}\end{author}
\begin{document}
\maketitle

\begin{abstract}\noindent
Let $E$ be an elliptic curve over  $\Q$. It is well known that the
ring of endomorphisms of $E_p$, the reduction of $E$ modulo a prime
$p$ of ordinary reduction, is an order of the quadratic imaginary
field $Q(\pi_p)$ generated by the Frobenius element $\pi_p$. When
the curve has complex multiplication (CM),  this is always a fixed
field as the prime varies. However, when the curve has no CM, very
little is known, not only about the order, but about the fields that
might appear as algebra of endomorphisms varying the prime. The ring
of endomorphisms is obviously related with the arithmetic of
$a_p^2-4p$, the discriminant of the characteristic polynomial of the
Frobenius element. In this paper, we are interested in the function
$\pi_{E,r,h}(x)$ counting the number of primes $p$ up to $x$ such
that $a_p^2-4p$ is square-free and in the congruence class $r$
modulo $h$. We give in this paper the precise asymptotic for
$\pi_{E,r,h}(x)$ when averaging over elliptic curves defined over
the rationals, and we discuss the relation of this result with the
Lang-Trotter conjecture, and with some other problems related to the
curve modulo $p$.
\end{abstract}

\section{Introduction and statement of results}

Let $E$ be an elliptic curve over  $\Q$ with conductor $N_E$.
For each prime $p$ of good reduction (i.e. $p \nmid
N_E$), $E$ reduces to an elliptic curve over the finite field
$\F_p$. The Frobenius endomorphism $(x,y) \mapsto (x^p, y^p)$ of
$E/\F_p$ has characteristic polynomial
$$x^2 - a_p  x + p=(x-\pi_p)(x-\overline{\pi}_p)$$
where $|a_p| \leq 2 \sqrt{p}$ by the Hasse bound. Let $\End(E/\F_p)$
be the ring of endomorphism of $E/\F_p$. If $p > 3$ is an ordinary
prime for $E$ (or equivalently, $a_p \neq 0$), then
$$\End(E/\F_p)\otimes\Q = \Q (\pi_p )$$
is completely determined by the Frobenius ring $\Z[\pi_p]$. The ring
of endormorphism $\End(E/\F_p)$ is more subtle, and the Frobenius
ring $\Z[\pi_p]$ can be a proper subset of $\End(E/\F_p)$. In fact,
it follows from Deuring's Theorem \cite{De} that any order $R$ such
that
$$\Z[\pi_p] \subseteq R \subseteq  \Q\left(\pi_p\right)$$ can
occur as the ring of endomorphism of some curve $E$ over $\F_p$ such
that the Frobenius of $E$ has characteristic polynomial $x^2 - a_p x
+ p$. It is then a natural question to ask whether $\Z[\pi_p] \simeq
\End(E/\F_p)$, or whether $\Z[\pi_p]$ is the maximal order of
$\Q(\pi_p)$ (the second question is of course a refinement of the
first one), and when this happens for a fixed curve $E$ when varying
the prime $p$.

For a fixed $E/\F_p$, let $\Delta_p$ be the discriminant of the
order $\End(E/\F_p)$. Then, the ``distance'' between the
endomorphism ring and the
 Frobenius ring is encoded by the unique positive integer $b_p$ defined by the relation $a_p^2 - 4 p
= b_p^2 \Delta_p$. The $b_p$'s where studied by Duke and Toth in
\cite{DuTo}, and by Cojocaru and Duke \cite{CoDu} who showed that
under the GRH
$$
\# \left\{ p  \leq x \;:\; p \nmid N_E \;\; \mbox{and} \;\;
b_p=1 \right\} \sim C_{\rm CD}(E) \frac{x}{\log{x}},
$$
for an explicit non-zero constant $C_{CD}(E)$. The condition $b_p=1$
is equivalent to the triviality of the Tate-Shafarevich group of $E$
over its function field $\F_p(E)$.

In general, properties of the Frobenius ring can be deduced from a
study of the arithmetic of $a_p^2-4p$ and, in particular, in
studying the square divisors of $a_p^2-4p$. In this paper, we study
the following question: given an elliptic curve $E$ over $\Q$, for
which primes $p$ is $a_p^2-4p$ square-free? For those primes, the
Frobenius ring $\Z[\pi_p]$ is the maximal order in $\Q(\pi_p)$. This
also provides a refinement of the question of  Cojocaru and Duke.
Since the result of  Cojocaru and Duke is itself a refinement of
Serre's result about the cyclicity of the group $E(\F_p)$, we have
the string of implications
$$a_p^2-4p \;\;\mbox{squarefree} \Longrightarrow b_p = 1 \Longrightarrow E(\F_p)
\;\;\mbox{cyclic}.$$

Can we then show that there is a positive proportion of $a_p^2-4p$
which are square-free, refining  Cojocaru and Duke's and Serre's
result? Maybe surprisingly, we believe that this is a much more
difficult question. If E has complex multiplication (CM), $a_p^2-4p$
square-free is equivalent to $p$ lying in some quadratic
progression. For example, let $E$ be the CM
elliptic curve $y^2=x^3-x$ with complex multiplication by $\Z[i]$.
Let $p$ be an ordinary prime, which is a prime congruent to 1 modulo
4. Since $E$ has rational 2-torsion, $a_p$ is even, and then 4
divides $a_p^2-4p$, and the natural question to consider is the square-freeness of
$(a_p^2-4p)/4$. Since $E$
has complex multiplication by $\Z[i]$, $a_p^2-4p = -4f^2$ for some
$f \in \Z$, and  $(a_p^2-4p)/4$ is square-free if and only if $f=1$
if and only if $p = (a_p/2)^2 + 1$.

Another reason to investigate square-freeness in the sequence $a_p^2
- 4p$ is that it might shed light on the following conjecture.
\begin{conj}[Lang-Trotter conjecture \cite{LaTr}]
\label{LTconjecture} Let $K$ be an imaginary quadratic number field,
and $E$ an elliptic curve over $\Q$ without complex multiplication.
Let
$$
\Pi_{E,K}(x) = \# \left\{ p \leq x \;:\; p \nmid N_E
\;\;\mbox{and}\;\; \Q(\pi_p) = K \right\}. $$
 Then as $x
\rightarrow \infty$
$$
\Pi_{E,K}(x) \sim C_{\rm LT}(E,K) \frac{\sqrt{x}}{\log{x}},$$ for
some explicit constant $C_{\rm LT}(E,K)$ depending on $E$ and $K$.
\end{conj}

Upper bounds for $\Pi_{E,K}(x)$  were investigated by Serre, by
Cojocaru, Fouvry and Murty \cite{CoFoMu}, and by Cojocaru and David
\cite{CoDa}. But there are no known lower bounds for $\Pi_{E,K}(x)$,
in particular there are no known examples of (non-CM) elliptic
curves such that $\Q(\pi_p) = K$ for infinitely many primes $p$. In
fact, we know much less than that, as there are no known examples of
(non-CM) elliptic curves such that $D_p$, the discriminant of
$\Q(\pi_p)$, lies in a fixed congruence class for infinitely many
primes $p$. Here again, one could give a lower bound for the number
of such primes by counting the number of primes $p$ such that
$a_p^2-4p$ is square-free and in the prescribed congruence class,
since in this case $a_p^2-4p=D_p$.

Let $h$ be a positive odd integer, and let $r$ be any integer such
that the greatest common divisor $(r,h)$ is square-free. Let
$$
\Pi_{E,r,h}^\be(x)= \#\{2<p\le x,\text{ prime
}\,:\,a_p^2-4p\in\Delta(r,h)\},
$$
where $\Delta(r,h)$ is the set of  square-free integers $n$ such
that $n \equiv r \mod h$. We remark that if $(r,h)$ is not
square-free, then $\Pi_{E,r,h}^\be(x) = 0.$  As for the restriction
to $h$ odd, it simplifies technical aspects of the proof in various
places, but it is not an essential restriction. Unfortunately, we
cannot give an asymptotic (or a lower bound) for
$\Pi_{E,r,h}^\be(x)$, but we can prove that the correct asymptotic
holds on average over elliptic curves.

In all the following,
$E(a,b)$ denotes the elliptic curve given by the equation
$y^2=x^3+ax+b$, with $4a^3 + 27 b^2 \neq 0.$

\begin{thm} \label{mainQ}
Let $h$ be a positive odd integer, and let $r$ be any integer such
that $(r,h)$ is square-free. Fix any $\varepsilon > 0$. Let $A, B$
be such that $AB
> x \log^{8}{x}$, $A, B > x^{\ep}$. Let $\mathcal{C}(A,B)$ be the set
of all elliptic curves $E(a,b)$ with integer coefficients $a,b$ such
that $|a| \leq A$ and $b \leq B$. Then, as $x \rightarrow \infty$,
$$
\frac{1}{\vert \mathcal{C}(A,B) \vert} \sum_{E(a,b) \in
\mathcal{C}(A,B)} \Pi_{E(a,b),r,h}^\be(x) = \mathfrak{C}
\frac{x}{\log{x}} + O \left( \frac{x}{\log^2{x}} \right),
$$
where $\mathfrak{C}$ is the positive constant
\begin{equation}\label{theconstant}
\mathfrak{C} = \frac{1}{3h} \prod_{{\ell \parallel h} \atop {\ell
\mid r}} \frac{\ell-1}{\ell} \prod_{{\ell \mid h}\atop{\ell \nmid
r}} \frac{\ell \left( \ell-1-\leg{r}{\ell} \right)}{(\ell-1)\left(
\ell-\leg{r}{\ell} \right)} \prod_{\ell \nmid h}
\frac{\ell^4-2\ell^2-\ell+1}{\ell^2(\ell^2-1)},
\end{equation}
where all products are taken over {\it odd} primes $\ell$ with the
specified conditions.
\end{thm}

The result above is similar in spirit, and by the techniques used to
prove it, to other average results as \cite{BaCoDa} for the
Koblitz's conjecture and \cite{FoMu} and \cite{DaPa} for the
Lang-Trotter conjecture. All those results rely on the fact that the
average over elliptic curves in ${\mathcal{C}}(A,B)$ can be
rewritten as an average over elliptic curves over $\F_p$ by
interchanging the summations, and this average can be reduced to an
average of class numbers by means of Deuring's Theorem. Hence,
Theorem \ref{mainQ} is equivalent to

\begin{thm} \label{mainFp}
Let $h$ be a positive odd integer, and let $r$ be any integer such
that $(r,h)$ is square-free. Let $$\Pi^\be(p)\ = \# \left\{
\mbox{$E$ over $\F_p$} \;:\; \mbox{$a_p^2-4p \in \Delta(r,h)$}
\right\}.$$ Then, as $x \rightarrow \infty$,
$$
\sum_{p \leq x} \Pi^\be(p)  = \frac{\mathfrak{C}}{3}
\frac{x^{3}}{\log{x}} + O \left( \frac{x^3}{\log^2{x}} \right),
$$
where $\mathfrak{C}$ is the constant of Theorem \ref{mainQ}.
\end{thm}


\

There is another sequence related to elliptic curves which was
investigated for square-freeness, namely the sequence $p+1-a_p =
\vert E(\F_p) \vert$ . Again, there should be a positive proportion
of square-free $p+1-a_p$, but this is still an open question, and
there are no known examples of non-CM elliptic curves over $\Q$ with
infinitely many square-free $p+1-a_p$. The case of CM curves is very
different. Let $E$ be an elliptic curve with CM by an order
${\mathcal O}$ in a quadratic imaginary field $K$, and let $p$ be a
prime of ordinary reduction. Then, $p+1-a_p$ is square-free if and
only if $\pi_p-1 = \mathfrak{p}_1 \dots \mathfrak{p}_s$ where the
$\mathfrak{p}_i$ are distinct split primes in ${\mathcal O}$, and
the number of such $\pi \in {\mathcal O}$ can be counted. When $E$
has CM by the maximal order of a quadratic imaginary field $K$, the
asymptotic for the number of primes $p$ such that $p+1-a_p$ is
square-free was obtained by Cojocaru \cite{Co}. The number of
square-free $p+1-a_p$ was also investigated in average over elliptic
curves $\F_p$  as $p$ tends to infinity by Gekeler \cite{Ge}. The
techniques that he uses are completely different from ours, and rely
on the work of Howe \cite{Ho} on the moduli spaces of elliptic
curves over $\F_p$ with a given group structure. Gekeler then proved
the analogue of our Theorem \ref{mainFp} for the sequence $p+1-a_p$
(his results are stated in a different way, but the two forms are
easily seen to be equivalent). As a corollary of his result
\cite[Proposition 4.2 (i)]{Ge}, one obtains, by going from the
average over $\F_p$ to the average over $\Q$ as indicated in Section
\ref{reduceaverage}, that the number of primes $p \leq x$ such that
$p+1-a_p$ is square-free follows the predicted asymptotic on average
over all elliptic curves. It is unclear if the techniques used by
Gekeler could be used to prove Theorem \ref{mainQ}, as $a_p^2-4p$ is
not easily related to the group structure of $E(\F_p)$.

The structure of our paper is as follows. In section
\ref{reduceaverage}, we reduce the average over elliptic curves over
$\Q$ to an average over elliptic curves over $\F_p$. There are
several ways to do that, and we use multiplicative characters as in
\cite{BaSh} and \cite{BaCoDa} to obtain the short average of Theorem
\ref{mainQ}. Section \ref{sectionbigsquares} deals with some easy
error terms of the asymptotic of Theorem \ref{mainQ}, while Section
\ref{sectionuseBDH} deals with the main error term by applying
 the Theorem of Barban, Davenport and Halberstam.
Section \ref{proofThm9} deals with the main term of the asymptotic,
and it is basically a computation of the constant
(\ref{theconstant}) of Theorem \ref{mainQ}. Finally, we explain in
Section \ref{sectionconjectures} what is the conjectural asymptotic
for $\Pi_{E,r,h}(x)$. The conjectural constant is very similar to
our average constant $\mathfrak{C}$, and Theorem \ref{mainQ} then
gives additional evidence for the conjecture.

\section{Reduction to a short average over finite fields}
\label{reduceaverage}

As we have mentioned in the introduction, the proof of Theorem
\ref{mainQ} relies in reducing the average of $\Pi_{E,r,h}^\be(x)$
to an average of class numbers by means of Deuring's Theorem. This
could be done in a straightforward manner by just  writing
\begin{eqnarray} \label{firstreduction}
\frac{1}{\vert \mathcal{C}(A,B) \vert} \ds\sum_{E \in
\mathcal{C}(A,B)} \Pi_{E,r,h}^\be(x) &=& \sum_{p\leq
x}\ds\sum_{(s,t) \in \mathcal{S}} \frac{1}{\vert \mathcal{C}(A,B)
\vert} \sum_{{|a|\leq A,\;|b|\leq B}\atop {{a\equiv s\pmod p}\atop
{b\equiv t\pmod p}}}1,
\end{eqnarray}
where $$\mathcal{S} = \mathcal{S}(p,r,h) = \left\{ (s_1,s_2) \in
\F_p^2 : a_p(E(s_1,s_2))^2-4p \in \Delta(r,h) \right\}.
$$

Observe that the number of terms in the middle sum is $\Pi^\be(p)$.
The innermost sum is simply
$$
\left(\frac{2A}p+O(1)\right)\left(\frac{2B}p+O(1)\right) \sim
\frac{4AB}{p^2}
$$
when $A, B$ are big enough with respect to $x$. Hence, we get  the asymptotic
$$
\frac{1}{\vert \mathcal{C}(A,B) \vert} \ds\sum_{E \in
\mathcal{C}(A,B)} \Pi_{E,r,h}^\be(x) \sim \sum_{p\leq
x}\frac{\Pi^\be(p)}{p^2}.
$$
However, this approach results a poor average, as we need to take
$AB > x^{2 + \epsilon}$ to get the asymptotic above.

Following \cite{BaSh} and \cite{BaCoDa}, one can obtain a
substantial improvement from a better use of the uniform
distribution of equivalent elliptic curves. Indeed, since two
elliptic curves $E(s_1,s_2)$ and $E(t_1,t_2)$ are isomorphic over
$\F_p$ if and only if $t_1=s_1u^4$ and $t_2=s_2u^6$ for some
$u\in\F_p^*$, there are $(p-1)/\#\Aut(E(s_1,s_2))$ elliptic curves
over $\F_p$ which are $\F_p$-isomorphic to a given elliptic curve
$E(s_1,s_2)$. Hence, we can write
\begin{eqnarray}\nonumber
&&\frac{1}{\vert \mathcal{C}(A,B) \vert} \ds\sum_{E \in
\mathcal{C}(A,B)} \Pi_{E,r,h}^\be(x) \\
\label{secondred} && \hspace{0.5cm} =\frac{1}{\vert \mathcal{C}(A,B)
\vert} \ds\sum_{p\leq x}\ds\sum_{(s,t) \in \mathcal{S}}
\frac{\#\Aut(E(s_1,s_2))}{p-1}\sum_{{|a|\leq A,\;|b|\leq B,\;\exists
1\leq u<p:}\atop {{a\equiv s_1 u^4\pmod p,}\atop {b\equiv s_2
u^6\pmod p}}}1,
\end{eqnarray}
where for each fixed $s_1, s_2 \in \F_p$, the innermost sum is over
all integers $|a| \leq A$, $|b| \leq B$ such that there exists $u
\in \F_p^*$ with $a \equiv s_1 u^4 \pmod p$ and $b \equiv s_2 u^6
\pmod p$. We now want to approximate the innermost sum by its main
term. This is the content of Lemma 4 of  \cite{BaCoDa}, and we will
use it to prove the short average of Theorem \ref{mainQ}. We include
it here for reading convenience.

\begin{lem}\cite[Lemma 4]{BaCoDa}
\label{nobaier:1} For any positive integer $k$,  we  have as $x
\rightarrow \infty$
\begin{eqnarray*}
&&\sum_{p\leq x}\frac1p\,\sum_{(s_1, s_2) \in (\F_p^*)^2}
\left|\sum_{{|a|\leq A,\;|b|\leq B,\;\exists 1\leq u<p:}\atop
{{a\equiv s_1u^4\pmod p}\atop {b\equiv tu^6\pmod
p}}}\kern-10pt1-\frac{2AB}{p}\right| \\  &&\ll
ABx^{1-\frac1{2k}}\log^{\frac k2-1}x + (A\sqrt B+B\sqrt
A)x^{1+\frac1{2k}}\log^{\frac k2-1}x + \sqrt{AB} x^{3/2}\log^2x.
\end{eqnarray*}
\end{lem}

Now, since
$$ \vert \mathcal{C}(A,B) \vert = 4AB + O(A+B),$$
and $\#\Aut(E(s,t))=2$ except when $p|a$ or $p|b$,
we can rewrite (\ref{secondred}) as
\begin{eqnarray}\label{lastEQ}
\frac{1}{AB} \sum_{p\leq x}\sum_{(s,t) \in
\mathcal{S}}\frac{1}{2(p-1)}\sum_{{|a|\leq A,\;|b|\leq B,\;\exists
1\leq u<p:}\atop {{a\equiv su^4\pmod p,}\atop {b\equiv tu^6\pmod
p}}} 1 +O\left( \log{\log{x}}+\frac{A+B}{AB}\frac {x}{\log x}\right)
\end{eqnarray}
using the trivial bound
$$
\frac{1}{\vert \mathcal{C}(A,B) \vert} \sum_{p \leq x} \sum_{{{|a|
\leq A} \atop {|b| \leq B}} \atop{ab \equiv 0 \pmod p}} 1 \ll
\log{\log{x}}+\frac{A+B}{AB}\frac {x}{\log x}.
$$
Then, using  Lemma \ref{nobaier:1} to replace the first term in
(\ref{lastEQ}), and observing that
the error term in (\ref{lastEQ}) is smaller than the error term of Lemma
\ref{nobaier:1}, we have
\begin{eqnarray}
&&\label{provethat}  \frac{1}{\displaystyle \vert \mathcal{C}(A,B)
\vert} \ds\sum_{E \in \mathcal{C}(A,B)} \Pi_{E,r,h}^\be(x) = \sum_{p
\leq x}
\frac{\Pi^\be(p)}{p(p-1)} \\
\nonumber && \kern-20pt+ O \left(
 x^{1-\frac1{2k}}\log^{\frac k2-1}x + (A^{-1/2}+
B^{-1/2}) x^{1+\frac1{2k}}\log^{\frac k2-1}x + (AB)^{-1/2}
x^{3/2}\log^2x \right). \end{eqnarray}

By choosing $A, B$ such that $A, B
> x^\varepsilon$ and $AB > x \log^{8}{x}$, and
$k$ large enough to have that $\varepsilon k
> 1$, we have that
\begin{eqnarray}
\label{provethat2} \frac{1}{\displaystyle \vert \mathcal{C}(A,B)
\vert} \ds\sum_{E \in \mathcal{C}(A,B)} \Pi_{E,r,h}^\be(x) &=&
\sum_{p \leq x} \frac{\Pi^\be(p)}{p(p-1)} + O \left(
\frac{x}{\log^2{x}} \right). \end{eqnarray}

The heart of the proof of Theorem \ref{mainQ} then consists in an
estimate for the main term in (\ref{provethat2}), and we will show
that
\begin{eqnarray}
\label{mainresultthm2} \sum_{p \leq x} \frac{\Pi^\be(p)}{p(p-1)} =
\mathfrak{C} \frac{x}{\log{x}} + O \left( \frac{x}{\log^2{x}}
\right),
\end{eqnarray}
which will complete the proof of Theorem \ref{mainQ}. With
(\ref{mainresultthm2}) and partial summation, we also get that
\begin{eqnarray}
\label{mainresultthm3} \sum_{p \leq x} \Pi^\be(p) =
\frac{\mathfrak{C}}{3} \frac{x^3}{\log{x}} + O \left(
\frac{x^3}{\log^2{x}} \right),
\end{eqnarray}
and the proof of Theorem \ref{mainFp} also follows from
(\ref{mainresultthm2}). We first apply Deuring's Theorem to rewrite
the sum in (\ref{provethat2}) as a sum of class numbers.

\begin{thm}(Deuring's Theorem \cite{De})\label{deuringthm}

\noindent Let $p > 3$ be a prime and let $t$ be an integer such that
$t^2 - 4p < 0$. Let $\mathcal{E}_r(p)$ be the set of
$\F_p$-isomorphism classes of elliptic curves defined over $\F_p$ with $a_p=t$.
Then,
$$
\ds\sum_{E \in \mathcal{E}_r(p)} \frac{1}{\# \Aut(E)} = H(t^2 - 4p),
$$
where $\Aut(E)$ is the automorphism group of $E$ and for any $D<0$,
$H(D)$ is the Kronecker class number
$$
H(D) := \sum_{{f^2 \mid D} \atop {\frac{D}{f^2} \equiv 0,1 (\mod 4)}}
\frac{h(D/f^2)}{w(D/f^2)}
$$
defined in terms of the class number $h(D/f^2)$ and number of
units $w(D/f^2)$ of $\Q(\sqrt{D/f^2})$.\\
Then, for any fixed $-2 \sqrt{p} \leq t \leq 2 \sqrt{p}$,
there are exactly $(p-1) H(t^2-4p)$ elliptic curves defined over
$\F_p$ with $a_p=t$.
\end{thm}
Using the previous theorem, and noting that for square-free $D$ we
have
$$
H(D) = \frac{h(D)}{w(D)},
$$
we can write
\begin{eqnarray} \label{forthm2}
\sum_{p \leq x} \frac{\Pi^\be(p)}{p(p-1)} = \sum_{{p \leq x, \; |t|
\leq 2 \sqrt{p}} \atop {t^2-4p\in\Delta(r,h)}}
\frac{h(t^2-4p)}{w(t^2-4p)p} = 2\sum_{p\le x}\sum_{{1\le t\le
2\sqrt{p}}\atop{t^2-4p\in \Delta(r,h)}}^{\kern14pt
odd}\frac{h(t^2-4p)}{w(t^2-4p)p}
\end{eqnarray}
since $t^2-4p$ is not square-free when $t$ is even.

\section{Terms divisible by big squares}
\label{sectionbigsquares}

 Let $\beta>0$ a parameter to be chosen later, and  let
\begin{eqnarray}\label{defKYR} K = [\log^{2\beta}x], \;\;Y = \frac{x}{K+1},
\;\;R=\log^{\beta}{x}.
\end{eqnarray}
We first split the sum over the primes up to $x$ into $K$ sums over
the intervals $[kY,kY+Y]$ of length $Y$, and write (\ref{forthm2})
as
\begin{eqnarray} \label{firstT}
\sum_{p \leq x} \frac{\Pi^\be(p)}{p(p-1)} = 2\sum_{1\le k\le
K}\sum_{{{kY< p\le kY+Y}}} \sum_{{1\le t\le 2\sqrt{p}} \atop {t^2-4p
\in \Delta(r,h)}} \sod \frac{h(t^2-4p)}{\omega(t^2-4p)p} + O\left(
\frac{x}{\log^{2\beta+1}{x}} \right),
\end{eqnarray}
where the error term is to account for the primes $p \leq Y$. Using
the Moebius function to detect squares, the main term of
({\ref{firstT}) is
\begin{eqnarray*}
2\sum_{1\le k\le K}\sum_{{{kY< p\le kY+Y}}}
\sum_{{1\le t\le 2\sqrt{p}}\atop {t^2-4p \equiv r \mod h}} \sod
\frac{h(t^2-4p)}{\omega(t^2-4p)p}  \sum_{d^2 \mid t^2-4p} \mu(d) &=&
T_1 + T_2 + T_3 ,
\end{eqnarray*}
where $T_1, T_2$ and $T_3$ will be sums over small, medium, and
large divisors $d$. In particular, for each $p$ we define the
intervals $I_1=[1,R]$, $I_2=(R,\sqrt Y]$, $I_3=(\sqrt Y,2\sqrt p]$
and then
\begin{eqnarray*}
T_i= 2\sum_{{1\le k\le K}\atop{kY< p\le kY+Y}} \sum_{d\in I_i}\sod\mu(d)\sum_{{{1\le
t\le 2\sqrt{p}}\atop{d^2|4p-t^2}}\atop {t^2-4p \equiv r \mod h}}
\sod\frac{h(t^2-4p)}{\omega(t^2-4p)p},
\end{eqnarray*}
for $i=1,2,3$. We show in this section that the terms corresponding to big squares
$d$, namely $T_2$ and $T_3$, are small and will be part of the error term.

For $T_3$ we use  the fact that there is a  unique representation of the numbers
$4p-t^2=m^2u$, where $u$ is
square-free, and that $d^2|4p-t^2$ if and only if $d|m$, to get  the trivial upper bound
$$
T_3\ll \log x\sum_{1\le k\le K}\sum_{kY< p\le
kY+Y}\frac{1}{\sqrt{p}}\sum_{\sqrt{Y}<m\le 2\sqrt p}\tau(m)|S_m(p)|,
$$
where $S_m(p)=\{1\le t\le 2\sqrt{p}\,:\, 4p-t^2=m^2u, u \text{
square-free}\}$, and $\tau(m)$ denotes the usual divisor function.
To get the above upper bound, and in other estimations to come, we
use the well-known bound $h(-d) \ll \sqrt{d} \log{d}$ for the class
number of a negative integer $-d$.
 Now, using the uniqueness of the representation, we can change
the sum in $m$ to a sum over square-free $u\le 4p/Y$ to obtain
\begin{eqnarray}\label{boundT3}
T_3&\ll& \frac{x^{\ep}}{Y}\sum_{1\le k\le K}\sum_{kY< p\le kY+Y}\sqrt{p}\ll
\frac{x^{3/2+\ep}}{Y},
\end{eqnarray}
where we have used that, for each square-free $u$ fixed, the number
of solutions to the equation $4p=t^2+m^2u$ is bounded above by the
number of ideal divisors of $4p$ in the ring of integers of
$\Q(\sqrt {-u})$, together with the trivial bound $\tau(m) \ll
m^\ep$. For $T_2$ we switch the sums to get
\begin{eqnarray}\nonumber
T_2&\ll& (\log x)\sum_{1\le k\le K}\frac{1}{\sqrt{kY}}\sum_{R<d\le
\sqrt{Y}}\sum_{{1\le t\le 2\sqrt{kY+Y}}}\sum_{{kY< p\le kY+Y}\atop{d^2|4p-t^2}}1\\
\label{boundT2}
&\ll& (\log x) Y\sum_{1\le k\le K}\sqrt{\frac{k+1}{k}}\sum_{R<d\le
\sqrt{Y}}\frac1{d^2}<
(\log x)\frac{x}{R}.
\end{eqnarray}

\section{An average of class numbers}
\label{sectionuseBDH}

The main term of $T$ comes from $T_1$. On each interval $kY < p \leq kY+Y$, we
approximate the sum over
$\{1\le t\le 2\sqrt{p}\}$ by $\{1\le t\le 2\sqrt{kY}\}$ with an error term of $RY
\log x$. Replacing (\ref{boundT2})
and (\ref{boundT3}) in (\ref{firstT}), we have
\begin{eqnarray*}
\sum_{p \leq x} \frac{\Pi^\be(p)}{p(p-1)}
&=&2\sum_{1\le k\le K}\sum_{{1\le t\le 2\sqrt{kY}}}\sod\sum_{d\le R}\sod\mu(d)
\sum_{{{kY< p\le kY+Y}\atop{d^2|t^2-4p}} \atop{t^2-4p \equiv r \mod h}}
\frac{h(t^2-4p)}{\omega(t^2-4p)p}\\
&& \;\;\;\;\;\;\;\;\; + O( x /\log^{\beta-1}x) \\
&=&\frac1\pi\sum_{1\le k\le K}\sum_{{1\le t\le 2\sqrt{kY}}}\sod\sum_{d\le
R}\sod\mu(d)\sum_{{{kY< p\le kY+Y}\atop{d^2|t^2-4p}} \atop{t^2-4p \equiv r \mod h}}
\frac{\sqrt{4p-t^2} \;L(1,\chi_{t^2-4p})}{p}\\
&& \;\;\;\;\;\;\;\;\; + O( x /\log^{\beta-1}x),
\end{eqnarray*}
using the
class number formula. Let
\begin{eqnarray}
\label{defU} U = {x^{1/2}} R^2.
\end{eqnarray}
Since
\begin{eqnarray*}
L(1,\chi_{t^2-4p})&=&\sum_{n\ge
1}\frac{\chi_{t^2-4p}(n)}{n}=\sum_{n\ge
1}\sod\frac{\chi_{t^2-4p}(n)}{n}+\frac12\sum_{{n\ge
1}}\frac{\chi_{t^2-4p}(2n)}{n}\\
&=&\left(1+\frac{\chi_{t^2-4p}(2)}2\right) \sum_{n\ge
1}\sod\frac{\chi_{t^2-4p}(n)}{n}+\frac14\sum_{{n\ge
1}}\frac{\chi_{t^2-4p}(n)}{n},
\end{eqnarray*}
using Polya-Vinogradov inequality to bound the tail of the
L-function, we obtain
\begin{eqnarray}
\nonumber
\sum_{p \leq x} \frac{\Pi^\be(p)}{p(p-1)} &=&\frac{2}{3\pi}\sum_{1\le k\le K}\sum_{{n\le U}\atop{1\le t\le
2\sqrt{kY}}}\sod\frac1n\sum_{d\le R}\sod\mu(d)
\sum_{{{kY< p\le kY+Y}\atop{d^2|t^2-4p}} \atop{t^2-4p \equiv r \mod h}}
\frac{\sqrt{4p-t^2}}{p}\chi_{t^2-4p}(n)\\
\label{collectT2T3}
&& \;\;\;\;\;\;\;\;\; + O( x /\log^{\beta-1}x),
\end{eqnarray}
since $\chi_{t^2-4p}(2)=-1$ for any $t,p$ in the previous conditions.
Using quadratic reciprocity, we can write $\chi$ as a character mod
$n$, and we use this to rewrite the main term of (\ref{collectT2T3}) as
\begin{equation*}
M
= \frac{2}{3\pi}\sum_{1\le k\le K}
\sum_{{{n\le U}\atop{1\le t\le
2\sqrt{kY}}}\atop{(t^2-r,h)=1}}
\sod\frac1n\sum_{{{\ap \tm{n}}\atop{(t^2-\ap,n)=1}}\atop {\ap \equiv r \mod (n,h)}}
\left(\frac \ap n\right)\sum_{{{d\le R} \atop {(d,nt)=1}}\atop{r \equiv 0 \mod
(d^2,h)}} \sod\mu(d)S(k,Y,n,t,\ap,d)
\end{equation*}
where
$$
S(k,Y,n,t,\ap,d)=\sum_{{{kY\le p\le kY+Y}\atop{p \equiv \nu(t,\ap,r) \mod
[nd^2,h]}}} \frac{\sqrt{4p-t^2}}{p},
$$
and $\nu(t,\ap,r)$ is the invertible residue modulo $[nd^2,h]$
solving the congruences
\begin{eqnarray*}
4p &\equiv t^2 - r& \mod h, \;\;\mbox{for $(t^2-r,h)=1$},\\
4p &\equiv t^2 - \alpha& \mod n, \;\;\mbox{for $(t^2-\alpha,n)=1$}, \\
4p &\equiv t^2& \mod d^2 \;\;\mbox{for $(t,d)=1$},
\end{eqnarray*}
under the conditions $(n,d)=1$, $r \equiv \alpha \mod(n,h)$ and $r \equiv 0 \mod
(d^2,h)$
insuring that the congruences are compatible.

In order to compute the asymptotic of $M$, we first change the weights and consider
\begin{equation*}
M' = \frac{2}{3\pi}\sum_{1\le k\le K} \sum_{{{n\le U}\atop{1\le t\le
2\sqrt{kY}}}\atop{(t^2-r,h)=1}} \sod\frac1n\sum_{{{\ap
\tm{n}}\atop{(t^2-\ap,n)=1}}\atop {\ap \equiv r \mod (n,h)}}
\left(\frac \ap n\right)\sum_{{{d\le R} \atop {(d,nt)=1}}\atop{r
\equiv 0 \mod (d^2,h)}} \sod\mu(d) S'(k,Y,n,t,\ap,d),
\end{equation*}
where
$$
S'(k,Y,n,t,\ap,d)=\sum_{{{kY\le p\le kY+Y}\atop{p
\equiv \nu(t,\ap,r) \mod [nd^2,h]}}} \frac{\sqrt{4kY-t^2}
\log{p}}{kY \log{kY}}.
$$
\begin{lem}
\label{changeweigth} Let $M$ and $M'$ be the two sums defined above
where $Y = x / (K+1)$ and $K=[\log^{2\beta}{x}]$ for some $\beta
\geq 1$. Then as $x \rightarrow \infty$
$$ M - M' \ll \frac{x}{\log^{\beta - 1}{x}}.$$
\end{lem}
\noindent{\bf Proof:}  It is enough to note that
\begin{multline}
\nonumber \kern-10ptS(k,Y,n,t,\ap,d)-S'(k,Y,n,t,\ap,d)=\kern-8pt\sum_{{{kY\le p\le
kY+Y}\atop{p
\equiv \nu(t,\ap,r) \mod [nd^2,h]}}} \left(\frac{\sqrt{4p-t^2}}{p}-\frac{\sqrt{4kY-t^2}
}{kY}\right)\\+\sum_{{{kY\le p\le kY+Y}\atop{p
\equiv \nu(t,\ap,r) \mod [nd^2,h]}}}\left(\frac{\sqrt{4kY-t^2}}{kY}-\frac{\sqrt{4kY-t^2}
\log{p}}{kY \log{kY}}\right)
\ll\frac{\sqrt Y}{knd^2}.
\end{multline}
Summing over $k,n,t,\ap,d$, with  trivial bounds for the remaining terms, gives the
result.
\qed

\

We then have to evaluate the sum $M'$ above. We first define some
notation.  For any positive integers $a,q$ such that $(a,q)=1$, let
\begin{eqnarray*}
\psi(X,Y;a,q) &=& \sum_{{X < p < X+Y} \atop {p \equiv a \tm q}}
\log{p} \\
E(X,Y;a,q) &=& \psi(X,Y;a,q) - \frac{Y}{\varphi(q)},
\end{eqnarray*}
where, as usual, $\varphi$ denotes the Euler function.

 Substracting  and adding $\displaystyle
\frac{Y}{\varphi([nd^2,h])}$ to the sum $M'$ above, we write it as
$M' = S_1 + S_2$ where
\begin{eqnarray}\label{firstS1}
S_1  = \frac{2Y}{3\pi}  \sum_{1\le k\le K} \sum_{{{n\le U} \atop
{{1\le t\le 2\sqrt{kY}}}}\atop{(t^2-r,h)=1}} \sod
\frac{\sqrt{4kY-t^2}}{kY \log{kY}} \frac{1}{n} \sum_{{{\ap
\tm{n}}\atop{(t^2-\ap,n)=1}}\atop{\alpha \equiv r \tm {(n,h)}}}
\leg{\ap}{n} \sum_{{{d\le R} \atop {(d,nt)=1}}\atop{r \equiv 0 \mod
(d^2,h)}} \frac{\mu(d)}{\varphi([nd^2,h])}
\end{eqnarray}
and
\begin{eqnarray}
\nonumber S_2 &\ll&  \sum_{1\le k\le K} \frac{1}{\sqrt{kY}}
 \sum_{{{1\le t\le 2\sqrt{kY}}\atop{(t^2-r, h)=1}}}\sod
\sum_{{n \leq U} \atop {{{d\le R}\atop{(d,nt)=1}}\atop{r \equiv 0
\tm {(d^2,h)}}}} \frac1n \sum_{{{\ap
\tm{n}}\atop{(t^2-\ap,n)=1}}\atop{\alpha \equiv r \tm {(n,h)}}}
\vert E(kY, Y; \nu(t, \alpha,r), [nd^2, h]) \vert
\\ \label{beforeCS}
&\ll& \sum_{{1\le k\le K}} \sum_{q \leq R^2 U h} \frac{R^3}{q}
\sum_{{a \mod q} \atop {(a,q)=1}} \vert E(kY, Y; a, q) \vert
\end{eqnarray}
since each $q \leq R^2 U h$ can be written as $q = [nd^2, h]$ in at
most $R$ ways, and for each fixed integer $t$, and for any $\alpha \not \equiv \alpha' \mod n$,
$\nu(t, \alpha,r) \not\equiv \nu(t, \alpha',r)$ modulo $q= [nd^2, h]$.
We now use:

\begin{thm}\label{BDH}(Barban-Davenport-Halberstam) Let $X, Y$ be positive integers
such that
$X+Y \leq x$. Then, for any $N > 0$, there exist $M > 0$ such that
$$
\sum_{q\le Q}\sum_{{a \mod q} \atop {(a,q)=1}}  \vert E(X, Y; a, q)
\vert^2 \ll \frac{x^2}{\log^N{x}}
$$
whenever $Q \leq x \log^{-M}{x}$.
\end{thm}

Using Cauchy-Schartz and Theorem \ref{BDH} on (\ref{beforeCS}), we
get that for any $A>0$
\begin{eqnarray} \nonumber
S_2 &\ll& R^3 \sum_{{1\le k\le
K}}  \left( \sum_{{{q \leq R^2
U h} \atop {a \mod q}} \atop {(a,q)=1}} \vert E(kY, Y; a, q) \vert^2
\right)^{1/2}  \left( \sum_{{{q \leq R^2 U h} \atop {a \mod q}}
\atop {(a,q)=1}} \frac{1}{q^2}
\right)^{1/2} \\ \label{ETforS2}
&\ll& x \log^{-A}{x}
\end{eqnarray}
since
$$ R^2 U h \le
\frac{x}{\log^M{x}}$$ for any $M>0$, and $x$ sufficiently large.

We now evaluate $S_1$. Consider the sum
\begin{equation} \label{defST}
S(T) =  \sum_{{1 \leq t \leq T}\atop{{(t^2-r,h)=1}}} \sod
\sum_{{{n\le U}}} \sod \frac{1}{n} \sum_{{{\ap
\tm{n}}\atop{(t^2-\ap,n)=1}}\atop{\alpha \equiv r \tm {(n,h)}}}
\leg{\ap}{n} \sum_{{{d\le R} \atop {(d,nt)=1}}\atop{r \equiv 0 \mod
(d^2,h)}} \frac{\mu(d)}{\varphi([nd^2,h])}.
\end{equation}

\begin{thm}
\label{mainresult} Let $R$, $U$ as before. Then, as $T,x \rightarrow
\infty$ we have for any $\ep > 0$
$$
S(T) =  \frac{3}{2} \mathfrak{C} T +O\left(\frac{T\log
R}{R}\right)+O\left(\frac1{T^{1-\ep}}\right).
$$
where $\mathfrak{C}$ is
the constant given in (\ref{theconstant}).
\end{thm}

Let us assume for a moment the previous theorem, which will be proved in the next
section, and let us show how the evaluation of $S_1$
in (\ref{firstS1}) follows from the asymptotic of $S(T)$. This will
complete the proof of Theorem \ref{mainQ}.

We define $F(t)$ such that
$$
S(T) = \sum_{{1 \leq t \leq T}\atop{{(t^2-r,h)=1}}} \sod F(t),
$$
where $S(T)$ is the sum defined in (\ref{defST}),
and consider $$ S_1(X) =
\sum_{{1 \leq t \leq 2\sqrt{X}}\atop{{(t^2-r,h)=1}}} \sod
\frac{\sqrt{4X-t^2}}{X
\log{X}} F(t).
$$
Using Theorem \ref{mainresult} and partial summation we get
\begin{eqnarray} \label{secondS1}
S_1(X) &=& \frac{\frac{3}{2} \mathfrak{C}}{X \log{X}} \int_{0}^{2 \sqrt{X}} t^2
(4X-t^2)^{-1/2} dt = \frac{ \frac{3}{2} \mathfrak{C}\pi}{\log X}+O(1/R).
\end{eqnarray}
Replacing (\ref{secondS1}) into (\ref{firstS1}), we get
\begin{eqnarray}\nonumber
S_1 &=& \frac{2Y}{3\pi}\sum_{1 \leq k \leq K} S_1(kY)
= \mathfrak{C}Y\sum_{1 \leq k
\leq K} \frac{1}{\log{kY}}+O(x/R) \\\label{above}
&=&
{\mathfrak{C}} \frac{x}{\log{x}}+O\left(\frac{x}{\log^2 x}\right).
\end{eqnarray}
Using (\ref{above}), and collecting the error terms from
(\ref{ETforS2}), Lemma \ref{changeweigth} and (\ref{collectT2T3})
with $\beta\ge 3$, this shows that
$$
\sum_{p \leq x} \frac{\Pi^\be(p)}{p(p-1)} =
{\mathfrak{C}} \frac{x}{\log{x}}+O\left(\frac{x}{\log^2 x}\right),$$
and it completes the proof of Theorem \ref{mainQ}.

\section{Proof of Theorem \ref{mainresult}}
\label{proofThm9}

Let
\begin{eqnarray*}
c_{t}(n) = \sum_{{{\alpha \mod n} \atop {(t^2-\alpha, n) =1}} \atop
{\alpha \equiv r \mod (n,h)}} \leg{\alpha}{n}
\end{eqnarray*}
be the sum over residues modulo $n$ that appears in (\ref{defST}).

\begin{lem} \label{aboutsmallc}
Let $r,h,t$ be integers  such that
$(t^2-r, h)=1$, and let $p$ be an odd prime. Then,
$c_{t}(n)$
is  a multiplicative function of $n$, with value at prime powers given by
\begin{itemize}{}{}
\item[1.] If $p \mid h$, then $c_t(p^\ell)p^{-\ell}= \displaystyle \frac{1}{(p ^\ell,h)}\left(\frac r{p^\ell}\right).$

\item[2.] If $p \nmid h$, then
\begin{eqnarray*}
c_t(p^\ell)p^{-\ell} = \left\{ \begin{array}{ll} 0 & \mbox{if $p \mid t$, $\ell$
odd;} \\
-1/p & \mbox{if $p \nmid t$, $\ell$ odd;} \\
1-1/p & \mbox{if $p \mid t$, $\ell$ even;} \\
1- 2/p & \mbox{if $p \nmid t$, $\ell$ even.} \end{array}
\right.
\end{eqnarray*}
\end{itemize}
\end{lem}
\noindent{\bf Proof:} This is a straightforward case by case computation, and we omit it.
\qed

\vspace{0.2cm}
First, for each fixed odd $t$ in (\ref{defST}), such
that $(t^2-r,h)=1$, we evaluate the sum
$$
\sum_{n \leq U}  \sod
 \frac{c_t(n)}{n} \sum_{{{d \leq R} \atop {(d,nt)=1}}\atop{r \equiv 0 \mod
(d^2,h)}}\sod \frac{\mu(d)}{\varphi([nd^2,h])}.
$$
For any $(d,n)=1$, we have that
\begin{eqnarray}\label{phi}
\varphi([nd^2,h]) =
\frac{\varphi(n)\varphi(d^2)\varphi(h)}{\varphi((n,h))
\varphi((d^2,h))}>\varphi(n)\varphi(d^2).
\end{eqnarray}
We also observe that
$$
\sum_{{{d> R} \atop {(d,nt)=1}}\atop{r \equiv 0 \mod (d^2,h)}}
\frac{1}{\varphi(d^2)}\ll\sum_{d> R} \frac{\log\log
d}{d^2}<\frac{\log R}{R},
$$
and
$$
\sum_{n> U}  \sod \frac{c_t(n)}{n\varphi(n)}\ll \sum_{n> U}  \sod
\frac{\log\log n}{n\prod_{p^{2\gamma+1}||n}p}\ll\sum_{n> U^{1/2}}
\sod \frac{2^{\omega(n)}\log\log n }{n^2}<\frac{1}{U^{1/2-\ep}}
$$
for any $\ep>0$. Then,
\begin{equation}\label{errorst}
\sum_{n\le U}  \sod
 \frac{c_t(n)}{n} \sum_{{{d\le R} \atop {(d,nt)=1}}\atop{r \equiv 0 \mod
(d^2,h)}} \sod \frac{\mu(d)}{\varphi([nd^2,h])}=
C_t+O\left(\frac{\log R}{R}\right),
\end{equation}
where
\begin{eqnarray}
\label{constant1} C_t =\sum_{n\ge 1}  \sod
 \frac{c_t(n)}{n} \sum_{{{d\ge 1} \atop {(d,nt)=1}}\atop{r \equiv 0 \mod
(d^2,h)}}\sod \frac{\mu(d)}{\varphi([nd^2,h])}=
\sum_{{{d\ge1}\atop{(d,t)=1}}\atop{(d^2,h)\mid r}}\sod \mu(d)
\sum_{{n\ge1}\atop{(n,d)=1}}\sod \frac{c_t(n)}{n \varphi([nd^2,h])}.
\end{eqnarray}
We now use the multiplicativity of the functions involved to
write the constant $C_t$ as an Euler product.
\begin{lem}
\label{eulerprod} Let $r,h,t$ be  integers such that $(t^2-r, h)=1$,
and $(h,r)$ is squarefree. Then,
$$\label{constant2} C_t =P(r,h)
\prod_{{p|t}\atop{p\nmid
h}}\left(1+\frac{2p^2+p-1}{p^4-p^3-2p^2-p+1}\right),
$$
where
$$
P(r,h)=\frac1{\varphi(h)}\prod_{p|h}\frac{p}{p-\left(\frac
rp\right)}\prod_{{p\parallel h}\atop{p \mid r}}
\left(1-\frac{1}{p}\right)\prod_{p\nmid
h}\frac{p^4-p^3-2p^2-p+1}{p(p-1)(p^2-1)},
$$
and all the products run over odd prime numbers.
\end{lem}
\noindent{\bf Proof:}
We first observe that
\begin{equation}\label{primefactor}
\frac{p^4-p^3-2p^2-p+1}{p(p-1)(p^2-1)}=1-\frac{p^2+2p-1}{p(p-1)(p^2-1)},
\end{equation}
so each of the products of  $C_t$ is convergent.
Although intricate,
the proof consists of a combination of straightforward computations
using Lemma \ref{aboutsmallc}. We always assume that every
prime appearing in the products below is odd. We first note that
$$
\sum_{{{d\ge 1} \atop {(d,nt)=1}}\atop{r \equiv 0 \mod (d^2,h)}}\sod
\frac{\mu(d)\varphi((d^2,h))}{\varphi(d^2)}=\prod_{{{p \nmid
nt}\atop{p
\parallel h}}\atop{p \mid r}}
\left(1-\frac{1}{p}\right) \prod_{p\nmid
nht}\left(1-\frac{1}{p(p-1)}\right)
$$
Hence, using (\ref{phi}), we get
$$
C_t=F(h,t)\sum_{n\ge 1}  \sod
\frac{c_t(n)\varphi((n,h))}{n\varphi(n)}\prod_{{{p|(r,n)}\atop{p
\parallel h}} \atop{p\nmid t}} \left(1-\frac{1}{p}\right)^{-1}
\prod_{{p|n}\atop{p\nmid ht}}\left(1-\frac{1}{p(p-1)}\right)^{-1},
$$
for
$$
F(h,t)=\frac{1}{\varphi(h)}\prod_{p\nmid
ht}\left(1-\frac{1}{p(p-1)}\right) \prod_{{{p
\parallel h} \atop {p \mid r}} \atop{p\nmid t}}\left(1-\frac{1}{p}\right).
$$
Again, writing the sum over $n$ as an Euler product, we get that
\begin{eqnarray} \label{cete}
C_t= F(h,t)\prod_{p}\left( 1 + \delta(p)\sum_{j\geq 1}
\frac{c_t(p^j)}{p^j} \frac{\varphi((p^j,h))}{\varphi(p^j)}\right),
\end{eqnarray}
where
\begin{eqnarray*}
\delta(p) = \left\{
\begin{array}{ll} \displaystyle \left(1 - \frac{1}{p(p-1)}\right)^{-1} &
\mbox{when $p \nmid ht$;} \\
  \\
 \displaystyle \left(1-\frac{1}{p}\right)^{-1} &
\mbox{when $p \parallel h, p \mid r, p \nmid t$;} \\
  \\
 \displaystyle 1 &
\mbox{otherwise.}
\end{array} \right.
\end{eqnarray*}
Finally, we use Lemma \ref{aboutsmallc} to note that when $p|h$, $
\displaystyle c_t(p^j)
\frac{\varphi((p^j,h))}{\varphi(p^j)}=\left(\frac rp\right)^j, $
independently of $t$. For primes $p\nmid h$, we have that
$\varphi((p^j,h))=1$. Hence, by splitting the inner sum in
(\ref{cete}) into odd and even terms, we get $$ \sum_{j\geq 1}
\frac{c_t(p^j)}{p^j}
\frac{\varphi((p^j,h))}{\varphi(p^j)}=\frac{p}{(p-1)(p^2-1)}\left(p\frac{c_t(p)}{p}+\frac{c_t(p^2)}{p^2}\right).
$$
The result now follows by considering the different cases in Lemma
\ref{aboutsmallc}.

 \qed

 \

We now proceed to prove Theorem \ref{mainresult}. In the following, all primes $p$ are odd and
all products
are restricted to odd primes.
Let us call $G$ the multiplicative function with value at a prime $p$ given by
$G(p)=\displaystyle\frac{2p^2+p-1}{p^4-p^3-2p^2-p+1}$.
Using (\ref{errorst}) and Lemma \ref{eulerprod}, we have that
\begin{eqnarray}\label{st}
\nonumber S(T)
\kern-8pt&=& \kern-8pt P(r,h)\sum_{{t \leq T} \atop
(t^2-r,h)=1}\sod\sum_{{d|t}\atop{(d,h)=1}}\mu^2(d)G(d)
+O\left(\frac{T\log R}{R}\right)\\
\kern-8pt&=& \kern-8pt P(r,h)\kern-4pt \sum_{{d\le
T}\atop{(d,h)=1}}\mu^2(d)G(d)\kern-8pt\sum_{{t \leq T/d} \atop (d^2t^2-r,h)=1}\sod
\kern-8pt 1+O\left(\frac{T\log R}{R}\right).
\end{eqnarray}
To compute the inner sum, we write
\begin{eqnarray}
\label{step1} \sum_{{t \leq T/d} \atop (d^2t^2-r,h)=1}\sod 1 &=&
\sum_{t \leq T/d}\sod\sum_{k| (d^2t^2-r,h)}\mu(k) = \sum_{k|h}
\mu(k) \sum_{{{t \leq T/d} \atop {t^2 \equiv rd^{-2} \mod k}} \atop
{t \equiv 1 \mod 2}} 1.
\end{eqnarray}
Because of the Moebius function, we can suppose that $k$ is
square-free. In that case, by the Chinese Remainder Theorem, the
number of solutions modulo $2k$ to the congruences $t^2 \equiv
rd^{-2} \mod k$ and $t \equiv 1 \mod 2$ is
$$ \prod_{p \mid k} \left(1 + \leg{r}{p}\right).$$
Using that in (\ref{step1}), we have that
\begin{eqnarray}\label{step2}
\sum_{{t\leq T/d} \atop (d^2t^2-r,h)=1}\sod 1 &=& \frac{T}{2d}
\left( \prod_{p \mid h}  1 - \frac{1 + \leg{r}{p}}{p} \right) +
O(\tau(h) 2^{\omega(h)}).\end{eqnarray}

  Now we observe that
\begin{equation}\label{second:sum}
\sum_{{d\le T}\atop{(d,h)=1}}\frac{\mu^2(d)G(d)}{d}=\prod_{p\nmid
h}\left(1+\frac{G(p)}{p}\right)+O\left(\frac1{T^{2-\ep}}\right).
\end{equation}

Using (\ref{step2}) and (\ref{second:sum}) in (\ref{st})  completes
the proof of Theorem \ref{mainresult} with
\begin{eqnarray*}
{\mathfrak C} &=& \frac{P(r,h)}{3} \prod_{p \mid h}
\left(1-\frac{1+\leg{r}{p}}{p}\right) \prod_{p \nmid h}\left( 1 +
\frac{G(p)}{p}\right) \\
&=& \frac{1}{2h} \prod_{{p \parallel h} \atop {p \mid r}}
\frac{p-1}{p} \prod_{{p \mid h}\atop{p \nmid r}}
\frac{p(p-1-\leg{r}{p})}{(p-1)(p-\leg{r}{p})} \prod_{p \nmid h}
\frac{p^4-2p^2-p+1}{p^2(p^2-1)}.
\end{eqnarray*}
\qed

\section{Conjectures and Constants}
\label{sectionconjectures}

We develop in this section a precise conjecture for the asymptotic
behavior of
$$
\Pi_{E,r,h}^\be(x)=\#\{2<p\le x,\text{ prime
}\,:\,a_p^2-4p \in\Delta(r,h)\},
$$
following the standard heuristics, as in \cite{LaTr}
or \cite{Ko} (for different but related questions). We want $a_p^2 -
4p$ to be square-free, i.e. not divisible by $\ell^2$, for any prime
$\ell$. If $\ell \mid h$, we also have to take into account the fact
that $a_p^2-4p \equiv r \mod h$.

Fix an odd prime $\ell$, and suppose for now that $\ell \nmid h$. To
measure the probability that $a_p^2-4p$ is not divisible by
$\ell^2$, we use the Galois extension $\Q(E[\ell^2])/ \Q$, where
$\Q(E[\ell^2])$ is the field obtained from $\Q$ by adjoining the
coordinates of all the $\ell^2$-torsion points of $E$. Since
$E[\ell^2]$, the group of $\ell^2$-torsion points of $E$, is isomorphic to
$\Z/\ell^2 \Z \times \Z / \ell^2 \Z$ as an
abstract group, fixing a basis of $E[\ell^2]$, we obtain an
injective group homomorphism
$$
\rho_{\ell^2} : \Gal(\Q(E[\ell^2])/ \Q) \hookrightarrow
\GL_2(\Z/\ell^2 \Z).
$$
Let $G(\ell^2)$ be the image of $\rho_{\ell^2}$ in $\GL_2(\Z/\ell^2
\Z).$ For any prime $p$ unramified in  $\Q(E[\ell^2])/ \Q$ (i.e. $p
\nmid \ell N_E$), the Frobenius automorphism $\sigma_p$ is sent
to a conjugacy class of matrices in $\GL_2(\Z/\ell^2 \Z)$ such that
for any $g$ in the conjugacy class
\begin{eqnarray*}
\tr(g) &\equiv& a_p \mod \ell^2\\
\det(g) &\equiv& p \mod \ell^2.
\end{eqnarray*}
Let $$\Omega(\ell^2) = \left\{ g \in G(\ell^2) :
\tr(g)^2 - 4 \det(g) \not\equiv 0 \mod \ell^2 \right\}.$$
We can then measure the probability that $\ell^2 \nmid a_p^2-4p$ by
the ratio
\begin{eqnarray*} P_1(\ell) &=& \frac{\vert \Omega(\ell^2) \vert}{\vert
G(\ell^2) \vert}.
\end{eqnarray*}
By Serre's Theorem \cite{Se}, $G(\ell^2) = \GL_2(\Z/\ell^2 \Z)$ for
all but finitely many primes $\ell$. Furthermore, there is an
integer $M_E$, depending on the elliptic curve $E$, such that the
Galois group of the extension obtained by adjoining the coordinates
of all torsion points of $E$ to $\Q$ is the full inverse image in
$GL_2(\hat{\Z})$ of $\Gal(\Q(E[M_E]/\Q) \subseteq \GL_2(\Z/M_E \Z)$,
where $\hat{\Z}$ is the ring of adeles of $\Z$. In particular, $\ell
\nmid M_E$ implies that $G(\ell^2) = \GL_2(\Z/\ell^2 \Z)$. The
integer $M_E$ is not uniquely defined, and if some integer $m$
satisfies the property above, so does any multiple of $m$. We remark
that $M_E$ is always even.

Then, for $\ell \nmid M_E$, we compute that
\begin{eqnarray*}
 P_1(\ell) &=& \frac{\vert \left\{ g \in \GL_2(\Z/\ell^2 \Z)  : \tr(g)^2 - 4 \det(g)
\not\equiv 0 \mod \ell^2 \right\} \vert}{\vert \GL_2(\Z/\ell^2 \Z) \vert}\\
&=& \frac{\ell^4-2\ell^2-\ell+1}{\ell^2(\ell^2-1)}.
\end{eqnarray*}

Suppose now that $\ell \nmid M_E$, but $\ell \mid h$. We write $h =
\ell^\alpha h'$ where $\alpha \geq 1$ and $(h', \ell)=1$. There are
two congruence conditions modulo powers of $\ell$, namely $a_p^2 -
4p \not\equiv 0 \mod \ell^2$ and $a_p^2 - 4 p \equiv r \mod
\ell^\alpha$. Let $\beta = \max(\alpha, 2)$. The probability that
$a_p^2-4p$ satisfies both congruences is then the ratio
$$
P_2(\ell) = \frac{\vert \left\{ g \in G(\ell^{\beta}) : \tr(g)^2 - 4
\det(g) \not\equiv 0 \mod \ell^2 \;\mbox{and}\; \tr(g)^2 - 4 \det(g)
\equiv r \mod \ell^\alpha \right\} \vert}{\vert G(\ell^\beta) \vert}.
$$
Again, Serre's theorem guarantees that $G(\ell^{\beta}) =
\GL_2(\Z/\ell^{\beta} \Z)$ since $\ell \nmid M_E$. For such an odd
prime $\ell$, one computes that the ratio above is
\begin{eqnarray*}
P_2(\ell) = \left\{ \begin{array}{ll} \displaystyle
\frac{\ell(\ell-1-\leg{r}{\ell})}{(\ell-1)(\ell - \leg{r}{\ell})
\ell^\alpha} & \mbox{when $\ell \nmid r$;}\\
\\ \displaystyle \frac{\ell-1}{ \ell^{\alpha+1}} & \mbox{when $\ell
\mid r$ and $\alpha=1$};\\
  \\
\displaystyle  \frac{1}{\ell^\alpha} & \mbox{when $\ell \parallel r$
and $\alpha \geq 2$;}\\
\\
\displaystyle  0 & \mbox{when $\ell^2 \mid r$
and $\alpha \geq 2$.}\\
\end{array} \right.
\end{eqnarray*}

The last case ($\ell^2 \mid r$ and $\alpha \geq 2$) happens only
when $\ell^2 \mid (r, h)$, which is impossible when $(r,h)$ is
square-free. If $(r,h)$ is not square-free, it is clear that
$\Pi_{E,r,h}^\be(x)$ is empty.


We now put all the local probabilities together. When $\ell_1,
\ell_2 \nmid M_E$, the probabilities modulo $\ell_1, \ell_2$ are
independent by Serre's Theorem. Using the probability $P_1(\ell)$
when $\ell \nmid hM_E$ and $P_2(\ell)$ when $\ell \nmid M_E, \ell
\mid h$, the conjectural constant is
\begin{eqnarray}
\label{conjecturalconstant} C_{SF}(E,r,h) = \frac{1}{h} P(M_E)
\prod_{{{\ell \nmid M_E}\atop{\ell
\parallel h}}\atop{\ell \mid r}} \frac{\ell-1}{\ell}
\prod_{{{\ell \nmid M_E}\atop{\ell \mid h}}\atop{\ell \nmid r}}
\frac{\ell(\ell-1-\leg{r}{\ell})}{(\ell-1)(\ell - \leg{r}{\ell})}
\prod_{{\ell \nmid M_E}\atop{\ell \nmid h}}
\frac{\ell^4-2\ell^2-\ell+1}{\ell^2(\ell^2-1)},
\end{eqnarray}
where $P(M_E)$ is the probability modulo $M_E$. In order to give an
explicit expression for $P(M_E)$, we need some notation. Let
$$
M_E = \prod_{\ell \mid M_E} \ell^{\gamma(\ell)}
\;\;\;\mbox{and}\;\;\; h = \prod_{\ell \mid h} \ell^{\alpha(\ell)} .
$$
If $\ell \mid M_E$, but $\ell \nmid h$, we define $\alpha(\ell)$ to
be $0$. For each prime $\ell \mid M_E$, let
\begin{eqnarray*}
\beta(\ell) =  \max{\left\{2, \gamma(\ell),\alpha(\ell)\right\}}.
\end{eqnarray*}
Finally, let
$$m = \prod_{\ell \mid M_E} \ell^{\beta(\ell)},$$
and let $G(m)$ be the image of the Galois representation
$$
\rho_{m} : \Gal(\Q(E[m])/ \Q) \hookrightarrow \GL_2(\Z/m \Z).
$$
Then,
$$
P(M_E)  = \frac{\vert \left\{ g \in G(m) : \tr(g)^2 - 4 \det(g)
\not\equiv 0 \mod \ell^{2} \;\mbox{and}\; \tr(g)^2 - 4 \det(g)
\equiv r \mod \ell^{\alpha(\ell)} \;\;\forall \ell \mid M_E \right\}
\vert}{\vert G(m) \vert}.
$$

\begin{conj} \label{squarefree1}
Let $h$ be a positive odd integer, and $r$ an integer such that
$(h,r)$ is square-free. Then as $x \rightarrow \infty$
$$
\Pi_{E(a,b),r,h}^\be(x) \sim C_{SF}(E,r,h) \frac{x}{\log{x}},
$$
where $C_{SF}(E,r,h)$ is the constant given by
(\ref{conjecturalconstant}).
\end{conj}

Then, Theorem \ref{mainQ} gives evidence for the asymptotic of
Conjecture \ref{squarefree1}, and also for the constant
$C_{SF}(E,r,h)$ appearing on the asymptotic.

A similar heuristic leads to the following conjecture for the
number of primes $p$ smaller than $x$ such that $p+1-a_p(E)$ is
square-free.
\begin{conj}
\label{squarefreeorders} Let $\Pi_E(x)$ be the number of primes of
good reduction such that $p+1-a_p(E)$ is square-free. Then as $x
\rightarrow \infty$
$$
\Pi_{E}(x) \sim C'_{SF}(E) \frac{x}{\log{x}},
$$
where $C'_{SF}(E)$ is the constant
$$
C'_{SF}(E) = \frac{|\Omega(M_E)|}{|G(M_E)|} \prod_{\ell \nmid M_E}  1
- \frac{\ell^3-\ell-1}{(\ell^2-1)\ell^2(\ell-1)}.
$$
\end{conj}

The work of Gekeler mentioned in the introduction then gives evidence for
Conjecture \ref{squarefreeorders}, and in particular for the constant
$C'_{SF}(E)$ appearing in the conjecture. We stress again that the
techniques used by Gekeler and our techniques are completely
different, but they both uncover the conjectural constants computed
from Galois representations in that section (and of course, neither
the work of Gekeler or the work of this present paper relies on
Galois representations, and counting elements of Galois groups).

\vspace*{1cm}

\noindent {\bf{Acknowledgements:}} Both authors would like to thank
Nathan Jones for useful discussions during the preparation of this
paper. This work was done during the visit of the second author at
the Centre de Recherches Math\'ematiques (CRM), and he thanks the
CRM for providing an excellent working environment.

\vspace{0.5cm}
\begin{small}
[Chantal David] Department of Mathematics and Statistics, Concordia
University, 1455 de Maisonneuve West, Montr\'eal, QC, Canada. {\it
Email:} cdavid@mathstat.concordia.ca \\

 [Jorge Jimenez Urroz] Universitat Polit\a`ecnica de Catalunya,
 Campus Nord, edifici C3, C. Jordi Girona, 1-3, 08034 Barcelona,
 Spain.
 {\it Email:} jjimenez@ma4.upc.edu
\end{small}

\end{document}